\documentclass[fleqn]{article}
\usepackage{graphicx}
\usepackage{amssymb,amsmath}
\usepackage{theorem}
\usepackage{hyperref}

\newtheorem{theorem}{Theorem}[section]
\newtheorem{question}[theorem]{Question}
\newtheorem{prop}[theorem]{Proposition}
\newtheorem{lemma}[theorem]{Lemma}

\newtheorem{corollary}[theorem]{Corollary}

\theorembodyfont{\rmfamily}

\newenvironment{proof}{\noindent{\bf Proof. }}{\hfill$\square$\medskip}

\begin{document}

\addtolength{\baselineskip}{3pt} \setlength{\oddsidemargin}{0.2in}

\def\R{{\mathbb R}}
\def\one{{\mathbf 1}}
\def\Q{{\mathbb Q}}
\def\Z{{\mathbb Z}}
\def\C{{\mathbb C}}
\def\N{{\mathbb N}}
\def\hom{{\rm hom}}
\def\Hom{{\rm Hom}}
\def\Inj{{\rm Inj}}
\def\Ind{{\rm Ind}}
\def\inj{{\rm inj}}
\def\ind{{\rm ind}}
\def\sur{{\rm sur}}
\def\PAG{{\rm PAG}}
\def\CUT{{\text{\rm CUT}}}
\def\eps{\varepsilon}
\long\def\killtext#1{}
\def\dd{{d_{\text{\rm set}}}}
\def\dl{{d_{\text{\rm left}}}}
\def\dr{{d_{\text{\rm right}}}}
\def\Ge{{\mathbf G}}

\def\ontop#1#2{\genfrac{}{}{0pt}{}{#1}{#2}}

\def\HH{{\cal H}}
\def\SS{{\cal S}}
\def\TT{{\cal T}}
\def\II{{\cal I}}
\def\FF{{\cal F}}
\def\GG{{\cal G}}
\def\AA{{\cal A}}
\def\MM{{\cal M}}
\def\EE{{\cal E}}
\def\LL{{\cal L}}
\def\PP{{\cal P}}
\def\CC{{\cal C}}
\def\SS{{\cal S}}
\def\QG{{\cal QG}}
\def\QGM{{\cal QGM}}
\def\CD{{\cal CD}}
\def\P{{\sf P}}
\def\E{{\sf E}}
\def\Var{{\sf Var}}
\def\T{^{\sf T}}

\def\tr{{\rm tr}}
\def\cost{\hbox{\rm cost}}
\def\val{\hbox{\rm val}}
\def\rk{{\rm rk}}
\def\diam{{\rm diam}}
\def\Ker{{\rm Ker}}

\title{\huge\bf Edge coloring models and reflection positivity\\[12mm]}

\author{{\sc Bal\'azs Szegedy}\\[8mm]}
%Microsoft Research \\
%One Microsoft Way\\
%Redmond, WA 98052}

%\date{Technical Report \\
%September 2004}

\maketitle

\tableofcontents

\begin{abstract}
Solving a conjecture of M. H. Freedman, L. Lov\'asz and A.
Schrijver we prove that a graph parameter is edge reflection
positive and multiplicative if and only if it can be represented
by an edge coloring model.
\end{abstract}

\section{Introduction}

The motivation of this paper comes from statistical physics as
well as from combinatorics and topology. The general setup in
statistical mechanics can be outlined as follows. Let $G$ be a
graph and let $\mathcal{C}$ be a finite set of "states" or
"colors". We think of $G$ as a crystal in which either the edges
or the vertices are regarded as "sites" which can have states from
$\mathcal{C}$. In the first case we speak about {\bf edge coloring
models} and in the second case about {\bf vertex coloring models}.
A configuration of the whole system is a function which associates
a state with each site.
 The states are interacting with each other at the vertices in edge
coloring models and along edges in vertex coloring models. A
weight is associated with each such interaction which is a real
(or complex) number depending on the interacting states (in vertex
coloring models there are additional weights associated the the
states). A concrete model is usually given by these numbers. The
partition function can be interpreted as a graph parameter which
is computed by summing the products of the weights over all
possible configurations of the system represented by $G$. It
proves to be useful to extend this graph parameter linearly to the
vector space of formal linear combinations of graphs. The elements
of this vector space are called quantum graphs. Quantum graphs
that can be obtained by gluing together a quantum graph with its
reflected version (using the distributive law) are called
reflection symmetric. However there are two different reasonable
definitions of gluing. In the first one we glue along unfinished
edges and in the second one along vertices. Correspondingly we get
the notions {\bf edge reflection symmetric} and {\bf vertex
reflection symmetric} quantum graphs. A graph parameter is called
{\bf edge reflection positive} (resp. {\bf vertex reflection
positive}) if it takes non negative values on edge-reflection
symmetric (resp. vertex reflection symmetric) quantum graphs. It
is a simple fact that the partition function in edge coloring
models is edge reflection positive and is vertex-reflection
positive in vertex coloring models. A surprising result proved by
M. H. Freedman, L. Lov\'asz and A. Schrijver (see \cite{FLS}) says
that vertex reflection positivity is almost enough to characterize
the partition functions of vertex coloring models. The extra
condition that they need is that the ranks of certain matrices
(which describe the gluing operation and are called connection
matrices) are growing at most exponentially. They conjectured that
similar characterization can be given for edge reflection positive
graph parameters. The main result of this paper (theorem
\ref{main}) is the proof of this conjecture in a strong version
where we replace the condition on the rank growth by a week and
natural condition namely that the graph parameter is
multiplicative for taking disjoint union of graphs.

The major difficulty of the proof is coming from a fact which is
interesting on its own right: In contrast with vertex coloring
models, partition functions of edge coloring models don't
determine the weights. There is an action of the orthogonal group
on different edge coloring models, which leaves the corresponding
partition function invariant. This phenomenon explains why it is
difficult to reconstruct an edge coloring model from its partition
function. In contrast with vertex coloring models we are searching
for an orbit of the orthogonal group rather than one specific
object. Our main tools to handle this difficulty are commutative
algebra and the theory of invariants of the orthogonal group.

A topological version of the above described reflection symmetry
and reflection positivity arises in topological quantum field
theory (see \cite{A} and \cite{FKNSW}) where the gluing operation
is defined on the formal linear space of manifolds with a fixed
boundary.

We should also emphasize that the subject has a close connection
to pure combinatorics. The partition function of a vertex coloring
model can be interpreted as the number of graph homomorphisms into
a fixed graph. This shows that the number of proper colorings and
many related important graph parameters are coming from vertex
coloring models (see \cite{FLS}). In many other cases where we
count certain structures in a graph (perfect matchings, fully
packed loop configurations etc...) it turns out that this number
is the value of the partition function of an edge coloring model.
The orthogonal invariance of edge coloring models generates
interesting equations between such numbers (A simple example is
shown in Chapter \ref{matchings}). Another peculiar fact, that we
show, is that vertex coloring models can be represented by complex
valued edge coloring models such that the values of the two
partition functions are identical. In some special cases the
representing edge coloring model is also real valued and in this
case the corresponding graph parameter is both vertex and edge
reflection positive. We show that the Ising model is such an
example. Finally We mention that a version of vertex coloring
models with an infinite number of states is worked out and
characterized in \cite{LSZ}. In such a model the states are
elements of a measure space on which the weights are given by a
measurable function. From the combinatorial point of view, these
vertex coloring models can be regarded as limits of sequences of
finite graphs and such objects are relevant to extremal
combinatorics. In Chapter \ref{ve} we point out that some of these
infinite models can be represented by edge coloring models with
finitely many states.

\section{Edge coloring models and reflection positivity}

\subsection{Circles and Quantum graphs}
Throughout the paper it will be convenient to extend the concept
of graphs by introducing edges that are not incident to any
vertex. We call such edges {\it circle} edges and picture them as
topological circles. Formally, a circle is an element of the edge
set which has no endpoints. Let $\mathcal{G}$ denote the set of
isomorphism classes of graphs, in which loops, circles and
multiple edges are allowed. We denote by $\emptyset\in
\mathcal{G}$ the empty graph whose vertex and edge sets are both
empty. If $G_1,G_2\in \mathcal{G}$ are two graphs then their
disjoint union $G_1\cup G_2$ is defined to be a graph whose vertex
set and edge set is the disjoint union of those of $G_1$ and
$G_2$. Every element of $\mathcal{G}$ is the disjoint union of an
ordinary (circle free) graph and a finite number of circles.

Let $R$ be an arbitrary commutative ring with $1$. An $R$-valued
graph parameter is a map $f:\mathcal{G}\rightarrow R$. We say that
$f$ is {\it multiplicative} if $f(G_1 \cup G_2)=f(G_1)f(G_2)$ for
any two graphs $G_1,G_2\in \mathcal{G}$ and $f(\emptyset)=1$.

Let $F$ be a field and let $\mathcal{Q}(F)$ denote the vectorspace
of finite $F-$linear combinations of elements of $\mathcal{G}$.
The elements of $\mathcal{Q}(F)$ are called {\it quantum graphs}.
The operation of taking disjoint union can be extended to quantum
graphs by using the distributive law. It is easy to see that
$\mathcal{Q}(F)$ becomes an $F$-algebra with $1$ if we introduce
disjoint union as multiplication. If the ring $R$ is an
$F$-algebra then any multiplicative graph parameter
$f:\mathcal{G}\rightarrow R$ extends uniquely to an algebra
homomorphism $f:\mathcal{Q}(F)\rightarrow R$. As a consequence we
have that the image of $\mathcal{Q}(F)$ is a subalgebra of $R$.

The usual setting throughout this paper is that $F$ is the field
of real numbers and $R$ is either the field of real numbers or a
polynomial ring over the reals. For this reason we use the
shorthand notation $\mathcal{Q}$ instead of
$\mathcal{Q}(\Bbb{R})$.

\subsection{Edge reflection positivity}\label{conn}

In this section we will need the notion of graphs with outgoing
(open) edges. An outgoing edge can be pictured as an edge which
goes out from the graph but is not finished. It can also happen
that such an edge goes out in both directions and so it is not
incident to any of the verteces. However the behavior of these
edges is different from circles because we want to maintain the
possibility of finishing them. To define this concept precisely we
need to introduce the set of "open ends" $O(G)$ of a graph $G$. A
graph $G$ with outgoing edges is a triple $(V(G),E(G),O(G))$ where
$(O(G)\cup V(G),E(G))$ is a graph in $\mathcal{G}$ with the
property that the degrees of open ends are exactly $1$.

We define $\mathcal{G}_k$ to be the set of all graphs with exactly
$k$ open ends which are labeled by the numbers $1,2,\dots,k$.
There is a natural operation
$$g:\mathcal{G}_k\times\mathcal{G}_k\rightarrow \mathcal{G}$$ which is
called ${\it gluing}$ and defined in the following way. Let $G_1$
and $G_2$ be two graphs in $\mathcal{G}_k$. Let us take the
disjoint union of them and identify their open ends which have the
same label. This way we obtain a graph in $\mathcal{G}$ in which
there are $k$ labeled vertices of degree $2$. Finally we eliminate
these vertices (and their two incident edges) by introducing a new
edge which connects their two neighbors directly. It is easy to
see that the resulting graph does not depend on the order in which
we eliminate the labeled vertices. Note also that the resulting
graph may contain circles even if the original two graphs did not
have any. This explains the importance of circles. The notation of
gluing is also defined for graphs with no outgoing edges, but in
this case gluing is the same as taking disjoint union.

Let $\mathcal{Q}_k$ denote the vectorspace of formal
$\Bbb{R}$-linear combinations of elements of $\mathcal{G}_k$. Now
the gluing operation extends uniquely to a symmetric bilinear form
$$g:\mathcal{Q}_k\times \mathcal{Q}_k\rightarrow \mathcal{Q}.$$
We say that a quantum graph $Q\in\mathcal{Q}$ is {\it edge
reflection symmetric} if $Q=g(H,H)$ for some quantum graph $H\in
\mathcal{Q}_k$ with $k\geq 0$. A graph parameter
$f:\mathcal{G}\rightarrow \Bbb{R}$ is called {\it edge reflection
positive} if its linear extension $f:\mathcal{Q}\rightarrow
\Bbb{R}$ takes non-negative values on all edge reflection
symmetric quantum graphs. In other words $f$ is edge reflection
positive if and only if the bilinear forms
$$f\circ g:\mathcal{Q}_k\times\mathcal{Q}_k\rightarrow \Bbb{R}$$
are positive semi-definite for all $k\geq 0$. One can write up the
matrices of these scalar products in the natural basis
$\mathcal{G}_k$ and obtain the so-called {\it connection matrices}
$M(k,f)$. These are infinite matrices whose rows and columns are
indexed by the elements of $\mathcal{G}_k$ and the entry in the
intersection of the row corresponding to $G_1$ and the column
corresponding to $G_2$ is $f(g(G_1,G_2))$.

\subsection{Edge coloring models and the characterization theorem}\label{edgmod}

Let $R$ be a commutative $\Bbb{R}$-algebra with $1$. (Usually
$R=\Bbb{R}$ or a polynomial ring over $\Bbb{R}$.) Let
$\mathcal{C}=\{c_1,c_2,\dots,c_d\}$ be a finite set of size $d$
whose elements will be referred as colors. An $R$-valued {\it edge
coloring model} is given by a function $t:\Bbb{N}^d\rightarrow R$
where $0$ is considered to be a natural number. For every edge
coloring model we are going to define an associated graph
parameter $t:\mathcal{G}\rightarrow R$. Let $v\in V(G)$ be a
vertex and let $\psi:E(G)\rightarrow \mathcal{C}$ be a coloring of
the edge set of a graph $G$. We denote by $v_{\psi}\in \Bbb{N}^d$
the vector whose $i$-th coordinate is the number of edges with
color $c_i$ incident to vertex $v$. It is important that loop
edges are counted twice. Now we define $t_{\psi}(G)$ by
$$t_{\psi}(G)=\prod_{v\in V(G)}t(v_{\psi})$$
and $t(G)$ by
$$t(G)=\sum_{\psi:E(G)\rightarrow \mathcal{C}}t_{\psi}(G).$$

It is clear that $t$ is a multiplicative graph parameter if we
define the empty product to be $1$ and moreover the value of $t$
on a single circle is the number of colors (which is $d$).

Let $k\geq 0$ be a natural number and let $G\in\mathcal{G}_k$ be a
graph with $k$ labeled outgoing edges. We say that a coloring
$\psi:E(G)\rightarrow\mathcal{C}$ is an extension of a coloring
$\chi:O(G)\rightarrow\mathcal{C}$ of the open ends if each open
end $o\in O(G)$ has the same color as the unique edge incident to
$o$. We denote this relation by $\psi>\chi$. For a coloring
$\chi:O(G)\rightarrow\mathcal{C}$ we introduce $t_{\chi}(G)$ by
$$t_{\chi}(G)=\sum_{\psi:E(G)\rightarrow\mathcal{C},~\psi>\chi}~~t_{\psi}(G)$$

Now let $G_1$ and $G_2$ be two graphs in $\mathcal{G}_k$. Since
the open ends in both $G_1$ and $G_2$ are labeled by numbers
$1,2,\dots,k$ we can say, by abusing the notation, that any
coloring $\chi:\{1,2,\dots,k\}\rightarrow\mathcal{C}$ is also a
coloring of $O(G_1)$ and $O(G_2)$. Assume that $\psi_1>\chi$ in
$G_1$ and $\psi_2>\chi$ in $G_2$ for the same coloring $\chi$.
Then there is a coloring $\psi=g(\psi_1,\psi_2)$ of the edges of
$G=g(G_1,G_2)\in\mathcal{G}$ which is obtained by gluing together
$\psi_1$ and $\psi_2$. This coloring has the property that
$v_{\psi}=v_{\psi_1}$ if $v\in V(G_1)$ and $v_{\psi}=v_{\psi_2}$
if $v\in V(G_2)$. It follows that
$$t_{\psi}(G)=t_{\psi_1}(G_1)t_{\psi_2}(G_2)$$
and that

\begin{equation}\label{szinszumm}
t(G)=\sum_{\chi:\{1,2,\dots,k\}\rightarrow\mathcal{C}}
t_{\chi}(G_1)t_{\chi}(G_2).
\end{equation}

 It is clear that the previous equality
also holds for $G_1,G_2\in\mathcal{Q}_k$ and
$G=g(G_1,G_2)\in\mathcal{Q}$ if we extend the invariants $t$ and
$t_{\chi}$ linearly to quantum graphs from $\mathcal{Q}$ and
$\mathcal{Q}_k$. As a consequence we get that real valued edge
coloring models give rise to edge reflection positive graph
parameters:

\begin{prop}
Let $t:\Bbb{N}^d\rightarrow\Bbb{R}$ be a real valued edge coloring
model. Then the graph parameter $t:\mathcal{G}\rightarrow \Bbb{R}$
is edge reflection positive.
\end{prop}

\begin{proof}
Let $k\geq 0$ ba a natural number and $Q=g(H,H)$ for some $H\in
\mathcal{G}_k$. Using equation (\ref{szinszumm}) we have that
$$t(Q)=\sum_{\chi:\{1,2,\dots,k\}\rightarrow\mathcal{C}}t_{\chi}(H)^2\geq
0.$$
\end{proof}

Our main theorem is the converse of the previous statement.

\begin{theorem}\label{main} Let $f:\mathcal{G}\rightarrow\Bbb{R}$ be an edge
reflection positive and multiplicative graph parameter. Then there
is an edge coloring model $t:\Bbb{N}^d\rightarrow\Bbb{R}$ such
that the corresponding graph parameter equals to $f$.
\end{theorem}

The subsequent chapters will lead to the proof of this theorem.

\subsection{Universal edge coloring models}

Let us fix a natural number $d$ and let us introduce algebraically
independent variables $x_v$ for each vector $v\in\Bbb{N}^d$. Let
$P_d$ denote the polynomial ring
$\Bbb{R}[\{x_v~|~v\in\Bbb{N}^d\}]$. The {\it universal edge
coloring model} $t_d$ corresponding to $d$ is a $P_d$ valued edge
coloring model which is given by the function $t_{d}(v)=x_v$. An
important property of these models is that real valued edge
coloring models $t$ with $d$ colors are in one to one
correspondence with homomorphisms $\varrho:P_d\rightarrow\Bbb{R}$
where the correspondence is given by the equation
$\varrho(x_v)=t(v)$. Note that if $t$ and $\varrho$ correspond to
each other then $t(Q)=\varrho(t_d(Q))$ for all $Q\in\mathcal{Q}$.

\medskip
Let us introduce $$I_d=\{t_d(Q)~|~Q\in\mathcal{Q}\}.$$ Since $t_d$
is multiplicative we have that $I_d$ is a subring of $P_d$. We
will prove later that $I_d$ is the set of all elements in $P_d$
which are invariant under a certain "natural" action of the
orthogonal group $O_d(\Bbb{R})$.

\subsection{Action of the orthogonal group on edge coloring
models}\label{aoe}

Let $d$ be a natural number and let $V$ be the vectorspace
consisting of the formal $\Bbb{R}$-linear combinations of the
colors $c_1,c_2,\dots,c_d$. We say that $V$ is the {\it color
space} and the elements of $V$ will be called {\it quantum
colors}. The space $V$ is endowed with an euclidean scalar product
for which $c_1,c_2,\dots,c_d$ is an orthonormal basis. Let us fix
an edge coloring model $t:\Bbb{N}^d\rightarrow R$. For every
natural number $n$ we define a symmetric $n$-linear form $l_n$ on
$V$ by
$$l_n(c_{i_1},c_{i_2},\dots,c_{i_n})=t(m_1,m_2,...,m_d)$$ where
$m_i$ denotes the number of occurrence of the color $c_i$ on the
list $c_{i_1},c_{i_2},\dots,c_{i_n}$.

Let $\alpha$ be an orthogonal transformation of $V$. We denote by
$u^{\alpha}$ the image of a vector $u\in V$ under the action of
$\alpha$. We define a new edge coloring model $t^{\alpha}$ by
$$t^{\alpha}(i_1,i_2,\dots,i_d)=l_n(c^{\alpha}_{j_1},c^{\alpha}_{j_2}
,\dots,c^{\alpha}_{j_n})$$ where $n=i_1+i_2+\dots+i_d$ and $d\geq
j_1,j_2,\dots,j_n\geq 1$ is an arbitrary sequence of integers such
that $|\{k~|~j_k=m\}|=i_m$. The goal of this section is to prove
the following.

\begin{prop}\label{ortin}
Let $G\in\mathcal{G}$ be an arbitrary graph. Then
$t(G)=t^{\alpha}(G)$.
\end{prop}

Let $G\in\mathcal{G}$ be a fixed circle free graph. A {\it half
edge} in $G$ is an ordered pair $(v,e)$ of a vertex $v$ and an
edge $e$ such that $v$ and $e$ are forming an incident pair. For
each edge $e\in E(G)$ we introduce two half edges $h(e,1)=(v_1,e)$
and $h(e,2)=(v_2,e)$ where $v_1$ and $v_2$ are the two endpoints
of $e$. In case $e$ is a loop, we think of $h(e,1)$ and $h(e,2)$
as different objects although the corresponding ordered pairs are
the same. We denote by $H(G)=\{h(e,i)~|~e\in E(G),~i\in \{1,2\}\}$
the set of half edges in $G$.

For each half edge $h(e,i)$, we introduce an isomorphic copy of
$V$ which we denote by $V_{e,i}$. In each such space $V_{e,i}$
there is a natural basis whose elements correspond to the colors
$c_1,c_2,\dots,c_d$. We denote the elements of this basis by
$c_{1,e,i},c_{2,e,i},\dots,c_{d,e,i}$. Let
$$W=\bigotimes_{e\in E(G),~i\in\{1,2\}}~V_{e,i}$$ be
the tensor product of all these spaces. For the edge coloring
model $t:\Bbb{N}^d\rightarrow R$ we define a linear form
$m:W\rightarrow R$ by
$$m\Bigl(\bigotimes_{e\in E(G)~,~i\in \{1,2\}}u_{e,i}\Bigr)=\prod_{v\in
V(G)}l_{d(v)}(u(v,1),u(v,2),\dots,u(v,d(v))$$ where $u_{e,i}\in
V_{e,i}$ are arbitrary elements, $d(v)$ is the degree of the
vertex $v$ and $u(v,1),u(v,2),\dots,u(v,d(v))$ is the list of
those $u_{e,i}$-s for which the half edge $h(e,i)$ is incident to
$v$. Since every half edge is incident to exactly one vertex we
have that the right hand side is multi linear in the vectors
$u_{e,i}$ and thus by the universal property of the tensor product
there is a unique $m$ which satisfies the equation.

Let us consider the spaces $W_e=V_{e,1}\otimes V_{e,2}$ associated
to the edges of $G$. A basis of $W_e$ is formed by the elements
$c_{i,e,1}\otimes c_{j,e,2}$ where $1\leq i,j\leq d$. Thus the
elements of $W_e$ can be represented as matrices whose rows and
columns are indexed by the elements of $\mathcal{C}$. Let
$J_e=\sum_{i=1}^d c_{i,e,1}\otimes c_{i,e,2}$ be the element of
$W_e$ which correspond to the identity matrix and let
$$J=\bigotimes_{e\in E(G)}J_e\in \bigotimes_{e\in E(G)}W_e=W.$$
We have that
$$J=\sum_{\psi:E(G)\rightarrow\{1,2,\dots,d\}}~\bigotimes_{e\in
E(G)~,~i\in\{1,2\}}c_{\psi(e),e,i}~.$$ Since the terms of this sum
correspond to the colorings of the edges in $G$ it follows that
$$t(G)=m(J).$$

Let $J_0$ denote the element $\sum_{i=1}^d c_i\otimes c_i$ in
$V\otimes V$.

\begin{lemma}\label{ortsumm}
If $b_1,b_2,\dots,b_d$ is an orthonormal basis in $V$ then
$J_0=\sum_{i=1}^d b_i\otimes b_i$.
\end{lemma}

\begin{proof}
Assume that $b_j=\sum_{i=1}^d a_{j,i}c_i$ for some real numbers
$a_{j,i}$. Then the matrix $A=\{a_{i,j}\}$ is orthogonal and thus
$AA^T=1$. Now $$\sum_{j=1}^d b_j\otimes b_j=\sum_{j=1}^d
\sum_{i=1}^d\sum_{k=1}^d a(j,i)a(j,k)c_i\otimes
c_k=\sum_{i=1}^d\sum_{k=1}^d \delta_{i,k}c_i\otimes c_k=J_0.$$
\end{proof}

Recall that $\alpha$ was an orthogonal transformation of $V$ and
thus $c_1^{\alpha},c_2^{\alpha},\dots,c_d^{\alpha}$ is an
orthonormal basis in $V$. By lemma \ref{ortsumm} we obtain that
$\sum_{i=1}^d c^{\alpha}_{i,e,1}\otimes c^{\alpha}_{i,e,2}=J_e$
and so

$$J=\sum_{\psi:E(G)\rightarrow\{1,2,\dots,d\}}~\bigotimes_{e\in
E(G)~,~i\in\{1,2\}}c^{\alpha}_{\psi(e),e,i}=J^{\alpha}$$ We obtain
that $$t(G)=m(J)=m(J^{\alpha})=t^{\alpha}(G)$$ for all circle free
graphs $G$.

Now let $H\in\mathcal{G}$ be an arbitrary graph which is the
disjoint union of a circle free graph $G$ and $n$ circles. The
equation $$t(H)=t(G)d^n=t^{\alpha}(G)d^n=t^{\alpha}(H)$$ completes
the proof of proposition \ref{ortin}.

\subsection{Action of the orthogonal group on the polynomial ring
$P_d$}\label{aop}

Recall that $P_d$ is the polynomial ring
$\Bbb{R}[\{x_v~|~v\in\Bbb{N}^d\}]$ and the universal edge coloring
model $t_d$ is given by the map $t_d:v\rightarrow x_v$. For a
fixed orthogonal transformation ${\alpha}$ of the color space $V$
we have a new edge coloring model
$t_d^{\alpha}:\Bbb{N}^d\rightarrow P_d$. Using that $P_d$ is a
free commutative $\Bbb{R}$-algebra with free generators
$\{x_v~|~v\in\Bbb{N}^d\}$ we get that the map
$$x_v\rightarrow t_d^{\alpha}(v)$$ extends to an algebra
endomorphism $\alpha:R_d\rightarrow R_d$. Since ${\alpha}^{-1}$
induces another endomorphism which is both right and left inverse
for ${\alpha}$ it turns out that ${\alpha}$ is an automorphism of
$R_d$. Proposition \ref{ortin} implies that

\begin{corollary}
The elements of the subring $I_d<R_d$ are invariant under the
action of ${\alpha}$ for all orthogonal transformations
${\alpha}$.
\end{corollary}

\medskip
We define the hight $h(x_v)$ of a variable $x_v\in P_d$ to be the
sum of the components of $v$. The hight of a monomial
$x_{v_1}x_{v_2}\dots x_{v_r}$ is defined to be the multiset
$\{h(x_{v_1}),h(x_{v_2}),\dots,h(x_{v_r})\}$. We denote by $W_S$
the linear subspace generated by all the monomials in $P_d$ of
hight $S$. It is clear that $P_d$ is the direct sum of the spaces
$W_S$ where $S$ runs over all possible finite multisets of the non
negative integers. We show that these subspaces are invariant
under the action of the orthogonal group. To see this let us fix
an element ${\alpha}\in O(V)$. Since $(x_{v_1}x_{v_2}\dots
x_{v_r})^{\alpha}=x_{v_1}^{\alpha}x_{v_2}^{\alpha}\dots
x_{v_r}^{\alpha}$ it is enough to prove that $x_v^{\alpha}$ is a
linear combination of variables of hight $n=h(x_v)$ for all $v$.
Let us represent $v$ by a multiset of colors
$\{c_{i_1},c_{i_2},\dots,c_{i_n}\}$. We have that
$$x_v^{\alpha}=t_d^{\alpha}(v)=l_n(c_{i_1}^{\alpha},c_{i_2}^{\alpha},\dots,c_{i_n}^{\alpha}).$$
By the multilinearity of $l_n$ the right hand side of the above
equation can be written as a linear combination of some monomials
of the form $l_n(c_{k_1},c_{k_2},\dots,c_{k_n})$ which are all
variables of hight $n$.

\section{Examples and questions}

\subsection{vertex coloring models as edge coloring
models}\label{ve}

A vertex coloring model (see \cite{FLS}) is given by a finite
weighted graph $H$ with real edge weights $\beta_H(i,j)$ and
positive vertex weights $\alpha_H(i)$. If $G$ is a simple graph
then the homomorphism function (or partition function) ${\rm
hom}(G,H)$ is defined by
$${\rm hom}(G,H)=\sum_{\phi:V(G)\to V(H)}~\prod_{v\in V(G)}\alpha_H(\phi(v))\prod_{uv\in
E(G)}\beta_H(\phi(u),\phi(v)).$$ We show that the graph parameter
${\rm hom}(G,H)$ can be represented by the partition function of
an edge coloring model where the number of colors is the rank of
the matrix of the edge weights in $H$. Let $B$ be the symmetric
matrix of the edge weights $\beta_H(i,j)$. From elementary linear
algebra we know that $$B=\lambda_1
u_1{u_1}^T+\lambda_2u_2{u_2}^T+\dots+\lambda_ru_r{u_r}^T$$ for
some real column vectors $u_i$ and numbers $\lambda_i\in\{1,-1\}$
where $r$ is the rank of $B$. Let $t$ be the edge coloring model
given by
$$t(s_1,s_2,\dots,s_r)=\sum_{j=1}^{|V(H)|}\alpha_H(j)\prod_{i=1}^r
(u_i(j)\sqrt{\lambda_i})^{s_i}.$$ Using that
$\beta_H(i,j)=\sum_{k=1}^{r}u_k(i)u_k(j)\lambda_k$ we get that
$t(G)={\rm hom}(G,H)$ for an arbitrary simple graph $G$. It is
worth mentioning that if $B$ is positive semidefinite then the
numbers $\lambda_i$ are all $1$ and the representing edge coloring
model is real valued. However there are cases when the edge
coloring model is real valued without this condition. A simple
example is the Ising model which can be represented by a weighted
graph on $2$ vertices with $\beta(1,1)=\beta(2,2)=a\geq 0$ ,
$\beta(1,2)=\beta(2,1)=b\geq 0$ and $\alpha(1)=\alpha(2)=1$. It is
easy to compute that the corresponding edge coloring model is
given by
$$t(s_1,s_2)=2\Bigl(\frac{a+b}{2}\Bigr)^{\frac{s_1}{2}}\Bigl(\frac{a-b}{2}\Bigr)^{\frac{s_2}{2}}~~{\rm if}~s_2~{\rm is~even}$$
$$t(s_1,s_2)=0~~{\rm if}~s_2~{\rm is~odd}.$$

It is an interesting phenomenon that the number of colors needed
to represent a vertex coloring model by an edge coloring model
depends only on the rank of the adjacency matrix of the weighted
graph. This leads to a family of infinite vertex coloring models
which are still representable by ordinary edge coloring models.
Let $w:[0,a]^2\rightarrow \Bbb{R}$ be a bounded symmetric
measurable function such that
$$w(x,y)=\sum_{i=1}^{r}\lambda_i f_i(x)f_i(y)$$ for some bounded
measurable functions $f_i$ and numbers $\lambda_i\in\{1,-1\}$.
Regarding the function $w$ as an infinite adjacency matrix one can
define an analogy of the homomorphism function by
$$t_w(G)=\int_{x_1,x_2,\dots,x_m}\prod_{(i,j)\in
E(G)~,~i<j}w(x_i,x_j)~dx_1~dx_2\dots dx_m$$ where $G$ is an
arbitrary graph with $|V(G)|=m$ such that the vertices of $G$ are
indexed by the numbers $\{1,2,\dots,m\}$. Note that $t_w$ is a
vertex reflection positive and multiplicative graph parameter. Let
us introduce the following edge coloring model
$$t(s_1,s_2,\dots,s_r)=\int_{x\in
[0,a]}\prod_{i=1}^{r}(f_i(x)\sqrt{\lambda_i})^{s_i}~dx.$$ It can
be calculated easily that $t_w(G)=t(G)$ and that if $\lambda_i=1$
for all $i$ then $t$ gives rise to a real valued edge coloring
model.

\subsection{Graph parameters from combinatorics}\label{matchings}

Many interesting graph parameters can be obtained from the
following special family of edge coloring models. Let $S$ be a
subset of $\Bbb{N}^r$ and let $t_S:\Bbb{N}^r\rightarrow\Bbb{R}$ be
the function with $t_S(v)=1$ if $s\in S$ and $t_S(v)=0$ if
$s\notin S$. Let $G$ be a simple graph. The next table lists a few
examples.

\begin{center}
\begin{tabular}{c|c}
  S & combinatorial meaning of $t_S(G)$\\
\hline
  $\{1\}\times\Bbb{N}$& number of perfect machings\\
  $\{2\}\times\Bbb{N}$& number of fully packed loop configurations\\
  $\{0,1\}\times\Bbb{N}$& number of matchins\\
  $\{0,2\}\times\Bbb{N}$& number of loop configurations\\
  $\{0,d\}\times\Bbb{N}$& number of $d$-regular subgraphs\\
  $\{0,1\}^d$& number of proper edge colorings with
$d$ colors\\
  $\{(b,c,d)|b+d\equiv c+d\equiv 0~(2)\}$& number of nowhere zero
4-flows\\
  $\{(2,0,0),(0,2,0),(0,0,1)\}\times\Bbb{N}$& permanent of the adjacency
matrix\\
\end{tabular}
\end{center}

Using the orthogonal invariance of partition functions one can
create peculiar equations. For example by rotating the firs
example on the list with $45$-degree we get the edge coloring
model
$$t(a,b)=\sqrt{2}^{~-(a+b)}(a-b)$$ whose partition function is
again the number of perfect matchings. In other words, the
partition function of the model $t(a,b)=a-b$ is $2^{|E(G)|}$ times
the number of perfect matchings in $G$.

\subsection{Open questions}

Let $f$ be an edge coloring model with $d$ colors and let $M(k,f)$
denote its $k$-th connection matrix (see chapter \ref{conn}). It
is not hard to see \cite{FLSP} that ${\rm rk}(M(k,f))\leq d^k$.
\begin{question} What are the possible sequences ${\rm
rk}(M(k,f))~,~k=1,2,3,\dots$?
\end{question}
The analogy of this question for vertex coloring models was
answered by L. Lov\'asz in $\cite{L}$.
\medskip

The next question is motivated by chapter \ref{ve}.
\begin{question} Which are those vertex coloring models whose
partition functions are edge reflection positive.
\end{question}
We know only two examples: The Ising model and the vertex coloring
models with positive semidefinite adjacency matrices.

\section{Proof of the characterization theorem}

Throughout this chapter we prove theorem \ref{main}. We assume
that $f$ is an edge-reflection positive and multiplicative graph
parameter. Recall that our goal is to construct an edge coloring
model $t:\Bbb{N}^d\rightarrow\Bbb{R}$ for some integer $d$ such
that the corresponding partition function equals to $f$.

\subsection{The value of a circle}

Let $k$ be a natural number and let $\mathcal{M}_{k}$ denote the
set of those graphs $G$ in $\mathcal{G}_{k}$ which are circle free
and $V(G)=\emptyset$. In particular the Edge set of $G$ is a
perfect matching on the $k$ open ends. It follows that if $k$ is
an odd number then $\mathcal{M}_k$ is empty. We denote by
$\mathcal{QM}_{k}$ the subspace generated by $\mathcal{M}_k$ in
$\mathcal{Q}_k$.

Assume that $k=2n$ for some natural number $n$ and let $A_k$
denote the subset of all elements $G$ of $\mathcal{M}_k$ with the
property that each edge of $G$ connects an open end with label
$\leq n$ and another open end with label $>n$. The elements of
$A_k$ can be represented by permutations of the set
$\{1,2,\dots,n\}$ in the way that a permutation $\pi$ correspond
to a matching $a_{\pi}\in A_k$ where the open end $i$ is connected
with $\pi(i)+n$ for all $1\leq i\leq n$. Now the definition of
gluing implies that $g(a_{\pi},a_{\varrho})$ is a graph which is
the disjoint union of $c(\pi\varrho^{-1})$ circles where
$c(\sigma)$ denotes the number of cycles in a permutation
$\sigma$. Recall  $d$ be the value of $f$ on a single circle.
Using the multiplicativity of $f$ we have that
$$f(g(a_{\pi},a_{\varrho}))=d^{~c(\pi\varrho^{-1})}.$$
Let $M_n$ be a matrix whose rows and columns are indexed by
permutations from the symmetric group $S_n$ and the entry in the
intersection if the row corresponding to $\pi$ and column
corresponding to $\varrho$ is $d^{~c(\pi\varrho^{-1})}$. The
assumption that $f$ is reflection positive implies that $M_n$ must
be a positive semidefinite matrix for every $n$. We will prove
that this is only possible if $d$ is a non negative integer.

The matrix $M_n$ is acting on the space of formal linear
combinations of the group elements of $S_n$ which is the group
algebra $\Bbb{R}[S_n]$. Let $$w=\sum_{\pi\in
S_n}\text{sgn}(\pi)\pi.$$ We have that
$$wM_k=\sum_{\pi,\varrho\in S_n}\text{sgn}(\pi)d^{~c(\pi\varrho^{-1})}\varrho=
\sum_{\pi,\varrho\in S_n
}\text{sgn}(\pi\varrho^{-1})d^{~c(\pi\varrho^{-1})}
\text{sgn}(\varrho)\varrho=w\Bigl(\sum_{\pi\in
S_n}\text{sgn}(\pi)d^{c(\pi)}\Bigr)$$ This means that $w$ is an
eigenvector of $M_n$ with eigenvalue $$\sum_{\pi\in
S_n}\text{sgn}(\pi)d^{c(\pi)}=d(d-1)(d-2)\dots(d-n+1).$$ The
positive semidefinitness of $M_n$ implies that
$d(d-1)\dots(d-n+1)$ must be a non negative number for every
natural number $n$ and so $d$ is a non negative integer.

\medskip
As the next lemma shows, the fact that the circle value is a non
negative integer is the first step towards the existence of an
edge coloring model representing $f$.

\begin{lemma}\label{korekviv}
Let $t$ be an arbitrary edge coloring model with $d$ colors. Then
$f(g(H,K))=t(g(H,K))$ for every pair $H,K\in \mathcal{QM}_k$.
\end{lemma}

\begin{proof}
The quantum graph $g(H,K)$ is the linear combination of graphs
consisting only of circles. The multiplicativity of $f$ shows that
the value of $t$ and $f$ must be the same on such a graph.
\end{proof}

Using the terminology of section \ref{edgmod}, we have that

\begin{lemma}\label{beingzero}
Let $t$ be any edge coloring model with $d$ colors and let $H\in
\mathcal{QM}_k$ then $f(g(H,H))=t(g(H,H))=0$ if and only if
$t_{\chi}(H)=0$ for all colorings
$\chi:\{1,2,\dots,k\}\rightarrow\mathcal{C}$.
\end{lemma}

\subsection{Lifting to the universal edge coloring model}\label{lifting}
Recall that $d$ is the circle value of $f$ and $t_d$ is the
universal edge coloring model with $d$ colors. In this section we
prove that

\begin{lemma}
If $Q\in\mathcal{Q}$ is an arbitrary quantum graph and $t_d(Q)=0$
then $f(Q)=0$.
\end{lemma}

\begin{proof}
For a graph $G\in\mathcal{G}$ we define its hight $h(G)$ to be the
multiset of the degrees of the vertices in $G$. From the
definition of the universal edge coloring model it follows that
$t_d(G)\in W_{h(G)}$. Every quantum graph $Q\in\mathcal{Q}$ can be
written in the form $\sum_{S}Q_S$ where $Q_S$ is a quantum graph
which is a linear combination of graphs of hight $S$. We have that
$t_d(Q)=\sum_{S}t_d(Q_S)$ and $t_d(Q_S)\in W_S$ for every multiset
$S$. It follow that $t_d(Q)=0$ implies that $t_d(Q_S)=0$ for all
multiset $S$.  Thus we can assume that $Q$ is homogeneous in the
sense that each graph component of $Q$ has the same hight $S$.

Let $G\in\mathcal{G}$ be a graph which is the disjoint union of a
circle free graph $H$ and $n$ circles. Both $f$ and $t_d$ vanish
on the quantum graph $G-d^nH$. This means that one can eliminate
circles in a quantum graph without changing the value of $f$ and
$t_d$ on it. Thus we can assume that $Q$ is a combination of
circle free graphs.

Assume that $S$ consists of $n$ numbers an their sum is $k$. Let
$G_S$ be a graph in $\mathcal{G}_k$ with the following properties:
\begin{enumerate}
\item $|V(G_S)|=n$
\item $|E(G_S)|=|O(G)|=k$
\item Each edge $e\in E(G)$ connects an open end with a vertex
\item The multiset of the degrees of the verteces is $S$.
\end{enumerate}

It is clear that $G_S$ is unique up to a relabeling of the open
ends. It is also clear that for every graph $G$ of hight $S$ there
is a matching $M\in \mathcal{M}_k$ such that $G=g(G_S,M)$. This
implies that our quantum graph $Q$ can be written in the form
$g(G_S,M)$ where $M$ is in $\mathcal{QM}_k$.

Let $\mathcal{P}_v\subseteq O(G)$ denote the set of those open
ends which are connected to the vertex $v\in V(G)$ in $G_S$ and
let $\mathcal{P}=\{\mathcal{P}_v|v\in V(G)\}$ be the partition
formed by these sets. We denote by $K\leq S_k$ the automorphism
group of $\mathcal{P}$. It is clear that $g(G_S,M)$ is isomorphic
to $g(G_S,M^{\sigma})$ for all ${\sigma}\in K$. Let
$$\hat{M}=\frac{1}{|K|}\sum_{\sigma\in K}M^{\sigma}.$$
We have that $$f(Q)=f(g(G_S,M))=f(g(G_S,\hat{M}))$$ and
$$0=t_d(Q)=t_d(g(G_S,M))=t_d(g(G_S,\hat{M}))$$.

Using equation (\ref{szinszumm}) from chapter \ref{edgmod} we get
that
$$0=t_d(Q)=\sum_{\chi:\{1,2,\dots,k\}\rightarrow\mathcal{C}}{t_d}_{\chi}(G_S){t_d}_{\chi}(\hat{M}).$$
The group $K$ is acting on the colorings
$\chi:\{1,2,\dots,k\}\rightarrow\mathcal{C}$. An orbit of this
action can be described as a multiset $X_{\chi}$ of multisets such
that the elements of $X_{\chi}$ are multisets of colors describing
the color distributions in different partition sets of
$\mathcal{P}$. The value of ${t_d}_{\chi}(G_S)$ depends only on
the orbit of $\chi$ because $G_S$ and $G_S^{\sigma}$ are
isomorphic for every ${\sigma}\in K$. Moreover ${t_d}_{\chi}(G_S)$
is a monomial of hight $S$ of the form $x_{v_1}x_{v_2}\dots
x_{v_n}$ where the vectors $v_i$ describe multisets of colors
which can be seen at different vertices and the list
$v_1,v_2,\dots,x_n$ describes $X_{\chi}$. On the other hand
${t_d}_{\chi}(\hat{M})$ is a real number which depend also only
the orbit of $\chi$ because $\hat{M}$ is $K$ symmetric. By using
the fact that different monomials are linearly independent over
$\Bbb{R}$ we obtain that ${t_d}_{\chi}(\hat{M})$ must be $0$ for
all colorings $\chi$. This implies by lemma \ref{beingzero} that
${t_d}(g(\hat{M},\hat{M}))=0$ and so by lemma \ref{korekviv} we
get that $f(g(\hat{M},\hat{M}))=0$.

Since $f(g(-,-))$ is a positive semidefinite form it follows that
$f(g(Y,\hat{M}))=0$ for all $Y\in\mathcal{G}_k$. In particular
$$f(Q)=f(g(G_S,\hat{M}))=0.$$
\end{proof}

\begin{corollary}\label{clifting}
There exists a homomorphism $\hat{f}:I_d\rightarrow\Bbb{R}$ such
that $f(Q)=\hat{f}(t_d(Q))$ for every quantum graph $Q\in
\mathcal{Q}$.
\end{corollary}

\begin{proof}
Recall that $\mathcal{Q}$ has an $\Bbb{R}$-algebra structure and
that $I_d$ is the image of $\mathcal{Q}$ under the algebra
homomorphism $t_d:\mathcal{Q}\rightarrow P_d$. On the other hand
$f:\mathcal{Q}\rightarrow\Bbb{R}$ is an algebra homomorphism
because $f$ is multiplicative. According to the main result of
this section we have that
$$\text{ker}(t_d)\subseteq\text{ker}(f)\subseteq\mathcal{Q}$$ and
this completes the proof.
\end{proof}

\subsection{Representing invariants of the orthogonal group with quantum
graphs}\label{riwq}

The main result of this section is the following.

\begin{lemma}\label{inrep}
If $p\in P_d$ is invariant under the action of the orthogonal
group $O_d(\Bbb{R})$ then there is a quantum graph
$Q\in\mathcal{Q}$ such that $t_d(Q)=p$.
\end{lemma}

Before we start proving lemma \ref{inrep} we describe a
construction which will be useful in this and in later sections.

\smallskip
\noindent{\bf tensor construction:}~~ Let $X$ be an arbitrary
finite set with a partition $\mathcal{Y}=\{Y_1,Y_2,\dots,Y_n\}$ on
its elements. Let $V_x$ be an isomorphic copy of the color-space
$V=\langle c_1,c_2,\dots,c_d\rangle_{\Bbb{R}}$ for each element
$x\in X$ and let
$$T(X,\mathcal{Y})=\bigotimes_{x\in X}V_x.$$
Let $l_n$ denote the symmetric $n$-linear form from chapter
$\ref{aoe}$ associated with the universal edge coloring model
$t_d$. For each partition set $Y_i$ we define a multilinear form
$\hat{m}_i$ by applying $l_{|Y_i|}$ for the spaces $\{V_x|x\in
Y_i\}$. The product $\prod_{i=1}^k \hat{m}_i$ defines a
multilinear form in the spaces $\{V_x|x\in X\}$. By factoring
$\hat{m}$ trough the tensor product $T(X,\mathcal{Y})$ we get an
$\Bbb{R}$-linear form $m:T(X,\mathcal{Y})\rightarrow P_d$. The
space $T(X,\mathcal{Y})$ admits a euclidean scalar product which
comes from the euclidean structure on $V$. An orthonormal basis
for this scalar product is formed by the different tensor products
of the color vectors. This basis will be called the {\it color
basis}. The orthogonal group $O_d(\Bbb{R})$ which preserves the
scalar product on $V$ is also acting on $T(X,\mathcal{Y})$ by
taking the tensor product of the actions on $V_x$. This action has
the property that $m(t^{\alpha})=m(t)^{\alpha}$ where $t\in
T(X,\mathcal{Y})$ and $\alpha\in O_d(\Bbb{R})$. Let
$S=\{|Y_1|,|Y_2|,\dots,|Y_n|\}$ be the multiset of the sizes of
the partition sets. If we substitute color vectors into the
multilinear form $\hat{m}$ we get all the monomial of hight $S$ in
$P_d$. It follows that $m$ maps $T(X,\mathcal{Y})$ to $W_s$
surjectively.
\medskip

\noindent{\bf Proof of lemma \ref{inrep}}~~ We know from chapter
\ref{aop} that $P_d$ is the direct sum (as a vectorspace) of the
spaces $W_S$ where each $W_S$ is invariant (as a subspace) under
the action of $O_d(\Bbb{R})$. It follows that if $p$ is an
invariant element of $O_d(\Bbb{R})$ then each $W_S$ component
$p_S$ of $p$ must be invariant too. Since $p$ is the sum of its
$W_S$ components it is enough to find quantum graphs $Q_S$ with
$t_d(Q_S)=p_S$ for each multiset $S$.

Let $S=\{s_1,s_2,\dots,s_n\}$ be a fixed multiset of natural
numbers and let $k=\sum_i s_i$. Let
$\mathcal{D}=\{D_1,D_2,\dots,D_n\}$ be a partition of the index
set $\{1,2,\dots,k\}$ such that $|D_i|=s_i$ for all $1\leq i\leq
n$. Let $W=T(\{1,2,\dots,k\},\mathcal{D})$ and let $W^0$ be the
kernel of the map $m$. We have that the space $W^0$ is invariant
under $O_d(\Bbb{R})$ and $W/W^0$ is isomorphic to $W_S$ in a way
that the induced action of $O_d(\Bbb{R})$ on $W/W^0$ commutes with
this isomorphism. By abusing the notation we identify $W_S$ with
$W/W^0$.

Let $p_1$ be a preimage of $p_S$ under the homomorphism
$W\rightarrow W/W^0$ and let
$$\bar{p}=\int_{\alpha\in O_d(\Bbb{R})}p_1^{\alpha}~d\nu$$
where $\nu$ is the normalized Haar measure on the orthogonal group
$O_d(\Bbb{R})$. Since $p_S$ is invariant in $W_S=W/W^0$ it follows
that $\bar{p}$ is also a preimage of $p_S$ under the map
$W\rightarrow W/W^0$. Furthermore we have that $\bar{p}$ is an
invariant of $O_d(\Bbb{R})$.

The first fundamental theorem of Weyl \cite{W} describes the space
of invariant elements in $W$ by determining a generating system
for it. The elements of this generating system correspond to
partitions of the set $\{1,2,\dots,k\}$ into two element subsets.
In particular if $k$ is an odd number then the only invariant is
the zero vector. Assume that $k$ is even and let
$\mathcal{E}=\{E_1,E_2,\dots,E_{k/2}\}$ be such a partition. Let
$\chi:\mathcal{E}\rightarrow\mathcal{C}$ be a coloring of the
partition sets. The function ${\chi}$ induces a coloring
$\hat{\chi}:\{1,2,\dots,k\}\rightarrow\mathcal{C}$ such that
$\hat{\chi}(j)=\chi(E_{i_j})$ for $1\leq j\leq k$ where $i_j$
denotes the number for which $j\in E_{i_j}$. We define $g_{\chi}$
to be the tensor product
$$\bigotimes_{i\in\{1,2,\dots,k\}}\hat{\chi}(i)$$ where
$\hat{\chi}(i)\in V_i$. The invariant which correspond to
$\mathcal{E}$ is
$$g=\sum_{\chi:\mathcal{E}\rightarrow\mathcal{C}}g_{\chi}.$$

We define a graph $G\in\mathcal{G}$ associated to the invariant
$g$. Let $V(G)=\{v_1,v_2,\dots,v_n\}$ and
$E(G)=\{e_1,e_2,\dots,e_{k/2}\}$. The edge $e_i$ connects the
vertices $v_{i_1}$ and $v_{i_2}$ where one element of $E_i$ is in
the partition set $D_{i_1}$ and the other element of $E_i$ is in
the partition set $D_{i_2}$. In other words, the vertices of $G$
correspond to the partition sets in $\mathcal{D}$, the edges
correspond to the partition sets in $\mathcal{E}$ and the edge
corresponding to $E_i$ is incident to the vertex corresponding to
$D_j$ if and only if $E_1\cap D_j\neq\emptyset$.

Now the spaces $V_i$ are in a one to one correspondence with the
half edges in $G$ and the form $m$ coincides with the one defined
in chapter \ref{aoe}. It follows that $m(g)=t_d(G)$.

Using Weyl's theorem we have that
$\bar{p}=\sum_{i=1}^r\lambda_ig_i$ for some real numbers
$\lambda_i$ and invariants $g_i$ where for each $g_i$ there is a
graph $G_i\in\mathcal{G}$ with $m(g_1)=t_d(G_i)$. It follows that
$$p_S=m(\bar{p})=\sum_{i=1}^r\lambda_i t_d(G_i)=t_d\Bigl(\sum_{i=1}^r\lambda_iG_i\Bigr).$$

\subsection{Projection to subalgebras of the matrix algebra}

Let $A$ be a subalgebra of the full matrix algebra
$\Bbb{M}_n(\Bbb{R})$ such that $A=\{M^T|M\in A\}$. The bilinear
function $(M,K)=\text{tr}(MK^T)$ defines a euclidean scalar
product on $\Bbb{M}_n(\Bbb{R})$. Let $\mathcal{P}_A$ denote the
orthogonal projection to $A$.

\begin{lemma}\label{projalg}
If $M$ is a symmetric positive semidefinite matrix then
$\mathcal{P}_A(M)=K^2$ for some symmetric matrix $K\in A$.
\end{lemma}

\begin{proof}
Since $A$ is invariant under transposing we have that
$\mathcal{P}_A(M)$ is a symmetric matrix in $A$. First we prove
that the eigenvalues of $\mathcal{P}_A(M)$ are all nonnegative.
Let $\lambda_1,\lambda_2,\dots,\lambda_k$ be the set of the
positive eigenvalues of $\mathcal{P}_A(M)$, let
$p(x)=(x-\lambda_1)(x-\lambda_2)\dots(x-\lambda_k)$ and let
$H=p(\mathcal{P}_A(M))$. Using that $\mathcal{P}_A$ is self
adjoint and that $H^T=H\in A$ have that
$$\text{tr}(\mathcal{P}_A(M)H^2)=(\mathcal{P}_A(M),H^2)=(M,\mathcal{P}_A(H^2))
=(M,H^2)=\text{tr}(MH^2)=\text{tr}(HMH)\geq 0.$$ Since
$\text{tr}(\mathcal{P}_A(M)H^2)$ is a positive linear combination
of the negative eigenvalues of $\mathcal{P}_A(M)$ it follows that
the eigenvalues of $\mathcal{P}_A(M)$ must be all nonnegative.

Let $g\in\Bbb{R}[x]$ be a polynomial such that
$g(\lambda_i)=\sqrt{\lambda_i}$ for $1\leq i\leq k$. Now
$K=g(\mathcal{P}_A(M))$ satisfy both $K^2=\mathcal{P}_A(M)$ and
$K\in A$.
\end{proof}

\subsection{Genaralized Brauer algebras}

Let $\mathcal{S}$ denote the set of finite multisets of the
positive integers. Let $\mu(S)$ denote the sum of the elements of
a multiset $S\in\mathcal{S}$. For each multiset $S\in\mathcal{S}$
we introduce a set $O(S)$ of size $\mu(S)$ and we define a
partition $P(S)$ on the elements of $O(S)$ such that the multiset
of the sizes of the partition sets in $P(S)$ is $S$. The algebra
$A_d$ consists of the formal linear combinations of triples
$$a(S_1,S_2,M)$$
where $S_1,S_2\in\mathcal{S}$ and $M$ is a perfect matching on the
set $O(S_1)\cup O(S_2)$. The product
$$a(S_1,S_2,M_1)a(S_3,S_4,M_2)$$
is defined to be $0$ if $S_2\neq S_3$. If $S_2=S_3$ then $M_1\cup
M_2$ is the edge set of a graph $G$ with node set $O(S_1)\cup
O(S_2)\cup O(S_4)$ such that nodes in $O(S_2)$ have degree $2$ and
nodes in $O(S_1)\cup O(S_4)$ have degree $1$. This means that $G$
is the union of node disjoint pathes and cycles. Replacing each
path by a single edge we get a matching $M_3$ on $O(S_1)\cup
O(S_4)$. Assume that the number of cycles in $G$ is $n$. The
product $a(S_1,S_2,M_1)a(S_2,S_4,M_2)$ is defined to be
$d^na(S_1,S_4,M_3)$.

We introduce the transpose map on $A_d$ as the unique linear
extension of the map $$a(S_1,S_2,M)^T=a(S_2,S_1,M).$$ Let
$A_d(S_1,S_2)$ denote the space spanned by the elements
$a(S_1,S_2,M)$ where $M$ runs through all perfect matchings of
$O(S_1)\cup O(S_2)$. We have that
$$A_d=\bigoplus_{S_1,S_2\in\mathcal{S}}A_d(S_1,S_2).$$
Let $$A_d(S)=\bigoplus_{S_1\in\mathcal{S}}A_d(S_1,S)$$ and
$$A_d(S)^T=\bigoplus_{S_1\in\mathcal{S}}A_d(S,S_1).$$
For an arbitrary basis element $a(S_1,S_2,M)$ we define
$\tau(a(S_1,S_2,M))\in\mathcal{G},~\tau_1(a(S_1,S_2,M))\in\mathcal{G}_{\mu(S_2)}$
and $\tau_2(a(S_1,S_2,M))\in\mathcal{G}_{\mu(S_1)}$ in the
following way. By identifying nodes in $O(S_1)\cup O(S_2)$
belonging to the same partition set of $P(S_1)\cup P(S_2)$ we get
$\tau(a(S_1,S_2,M))$. By identifying nodes in $O(S_1)$ (resp
$O(S_2)$) belonging to the same partition set of $P(S_1)$ (resp.
$P(S_2)$) and defining $O(S_2)$ (resp. $O(S_1)$) to be the set of
open edges we get $\tau_1(a(S_1,S_2,M))$ (resp.
$\tau_2(a(S_1,S_2,M))$. The map $\tau$ extends linearly to a map
$\tau:A_d\rightarrow \mathcal{Q}$ and the maps $\tau_1$,~$\tau_2$
extend to maps
$$\tau_1:A_d(S)\rightarrow\mathcal{Q}_{\mu(S)}~,~\tau_2:A_d(S)^T\rightarrow\mathcal{Q}_{\mu(S)}.$$

\begin{lemma}\label{posit}
If $b\in A_d$ then $f(\tau(bb^T))\geq 0$.
\end{lemma}

\begin{proof}
We use that $$b=\sum_{S\in\mathcal{S},\mu(S)\leq m}b_S$$ where
$b_S\in A_d(S)$ for all $S$ and $m$ is a large-enough natural
number. Since $$bb^T=\sum_{S\in\mathcal{S},\mu(S)\leq m}b_Sb_S^T$$
it suffices to show that $f(\tau(b_Sb_S^T))\geq 0$ for all $S$.
Using that $\tau_1(b_S)=\tau_2(b_S^T)$ we obtain that
$Q=g(\tau_1(b_S),\tau_2(b_S^T))$ is a reflection symmetric quantum
graph. The graph $\tau(b_Sb_S^T)$ can be obtained from $Q$ by a
process where in each step we delete a circle from a graph
component and multiply it by $d$. Using that the circle value of
$f$ is $d$, and that $f$ is multiplicative and edge reflection
positive we have that
$$f(\tau(b_Sb_S^T))=f(Q)\geq 0.$$
\end{proof}

Now we describe a matrix representation of the algebra $A_d$ which
will be of crucial importance in the next section. Let us
introduce the notation $\Bbb{M}(X,Y)$ for the space of real
matrices whose rows are indexed by the set $X$ and whose columns
are indexed by the set $Y$ where $X$ and $Y$ are finite sets.
Recall that $\mathcal{C}=\{c_1,c_2,\dots,c_d\}$ is a set with $d$
colors. Let $S_1,S_2\in\mathcal{S}$ be two multisets and let $M$
be a perfect matching on $O(S_1)\cup O(S_2)$. We say that a
coloring of $O(S_1)\cup O(S_2)$ is compatible with $M$ if the two
endpoints of each matching edge have the same color. First we
represent $a(S_1,S_2,M)$ by a matrix whose rows are indexed by
colorings $O(S_1)\rightarrow\mathcal{C}$ and whose columns are
indexed by colorings $O(S_2)\rightarrow\mathcal{C}$. The entry in
the intersection of the row $\chi$ and column $\psi$ is $1$ if the
coloring $\chi\times\psi:O(S_1)\cup O(S_2)\rightarrow\mathcal{C}$
is compatible with $M$ and is $0$ otherwise. By extending this
representation linearly to $A_d(S_1,S_2)$ we obtain a map
$$\omega:A_d(S_1,S_2)\rightarrow\Bbb{M}(\mathcal{C}^{O(S_1)}~,~\mathcal{C}^{O(S_2)}).$$
The reader can check easily that the map $\omega$ satisfies the
identity
$$\omega(a(S_1,S_2,M_1)a(S_2,S_3,M_2))=\omega(a(S_1,S_2,M_1))\omega(a(S_2,S_3,M_2)).$$
It follows that $\omega(b_1b_2)=\omega(b_1)\omega(b_2)$ if $b_1\in
A_d(S_1,S_2)$ and $b_2\in A_d(S_2,S_3)$. Let
$$\hat{A}_d(S_1,S_2)=\Bbb{M}(\mathcal{C}^{O(S_1)}~,~\mathcal{C}^{O(S_2)})$$
and let
$$\hat{A}_d=\bigoplus_{S_1,S_2\in\mathcal{S}}\hat{A}_d(S_1,S_2).$$
The space $\hat{A}_d$ is endowed with a natural algebra structure
in the following way. Assume that
$$b_1\in\hat{A}_d(S_1,S_2)~,
~b_2\in\hat{A}_d(S_3,S_4).$$ If $S_2=S_3$ then $b_1b_2$ is the
usual matrix product and if $S_2\neq S_3$ then $b_1b_2$ is defined
to be $0$. This multiplication rule defines a multiplication on
the whole space $\hat{A}_d$. It is clear that the map $\omega$
extends to an algebra homomorphism
$\omega:A_d\rightarrow\hat{A}_d$. Let
$$A_{d,r}=\bigoplus_{S_1,S_2\in\mathcal{S},~\mu(S_1),\mu(S_2)\leq
r}A_d(S_1,S_2)$$ and let
$$\hat{A}_{d,r}=\bigoplus_{S_1,S_2\in\mathcal{S},~\mu(S_1),\mu(S_2)\leq r}\hat{A}_d(S_1,S_2)$$
The space $A_{d,r}$ is a subalgebra of $A_d$ and the space
$\hat{A}_{d,r}$ is a subalgebra of $\hat{A}_d$. Moreover, $\omega$
maps $A_{d,r}$ into $\hat{A}_{d,r}$. It is easy to see that
$\hat{A}_{d,m}$ is the full matrix algebra
$$\Bbb{M}\Bigl(\bigcup_{S\in\mathcal{S},~\mu(S)\leq r}\mathcal{C}^{O(S)}~,~\bigcup_{S\in\mathcal{S},~\mu(S)\leq
r}\mathcal{C}^{O(S)}\Bigr)$$ and that $\omega(b^T)=\omega(b)^T$
for all $b\in A_{d,r}$. This implies in particular that
$\omega(A_{d,r})$ is a subalgebra of the matrix algebra
$\hat{A}_{d,r}$ which is closed under taking transpose. Let us
define the euclidean scalar product $(b_1,b_2)=\text{tr}(b_1b_2)$
on $\hat{A}_{d,r}$. The spaces $\hat{A}_d(S_1,S_2)$ are orthogonal
to each other in $\hat{A}_{d,r}$ for different pairs $(S_1,S_2)$.
Since $\omega(A_d(S_1,S_2))$ is contained in $\hat{A}_d(S_1,S_2)$
we have that
$$\omega(A_{d,r})=\bigoplus_{S_1,S_2\in\mathcal{S},~
\mu(S_1),\mu(S_2)\leq r}\omega(A_d(S_1,S_2))$$ where all the
direct summands are orthogonal to each other. Let
$\mathcal{P}_{d,r}$ denote the orthogonal projection of
$\hat{A}_{d,r}$ to $\omega(A_{d,r})$ and let
$\mathcal{P}_{d,r}(S_1,S_2)$ denote the orthogonal projection of
$\hat{A}_{d,r}$ to $\omega(A_d(S_1,S_2))$ where
$\mu(S_1),\mu(S_2)\leq r$. The above properties imply that the
restriction of $\mathcal{P}_{d,r}$ to the space
$\hat{A}_d(S_1,S_2)$ is $\mathcal{P}_{d,r}(S_1,S_2)$.

Now we use the tensor construction from chapter \ref{riwq}. Let us
observe that the elements of the color basis in the space
$T(O(S_1)\cup O(S_2),P(S_1)\cup P(S_2))$ correspond to colorings
$O(S_1)\cup O(S_2)\rightarrow\mathcal{C}$ and the elementary
matrices in $\hat{A}_d(S_1,S_2)$ also correspond to such
colorings. This gives an isometry between the euclidean spaces
$T(O(S_1)\cup O(S_2),P(S_1)\cup P(S_2))$ and $\hat{A}_d(S_1,S_2)$.
By abusing the notation, we will identify the these spaces. Now
the tensor construction defines maps
$m:\hat{A}_d(S_1,S_2)\rightarrow P_d $ for all
$S_1,S_2\in\mathcal{S}$. These maps have a unique common linear
extension $m:\hat{A}_d\rightarrow P_d$.

Let us observe that $$T(O(S_1)\cup O(S_2),P(S_1)\cup
P(S_2))=T(O(S_1),P(S_1))\otimes T(O(S_2),P(S_2))$$ and that for
$v_1\in T(O(S_1),P(S_2)),~v_2\in T(O(S_2),P(S_2))$ we have that
$$m(v_1)m(v_2)=m(v_1\otimes v_2).$$ Let
$$\hat{B}_{d,m}=\bigoplus_{S\in\mathcal{S},~\mu(S)\leq
m}T(O(S),P(S)).$$ Is is clear that
$$\hat{B}_{d,m}\otimes\hat{B}_{d,m}=\hat{A}_{d,m}$$ and that
$$m(v_1\otimes v_2)=m(v_1)m(v_2)$$ for an arbitrary pair
$v_1,v_2\in\hat{B}_{d,m}$.

% Let
%$S_1,S_2\in\mathcal{S}$ be two multisets and let us introduce
%isomorphic copies of the color-space $V$ for each element in
%$O(S_1)\cup O(S_2)$. We denote by $T(S_1,S_2)$ the tensor product
%of these spaces. Note that the space $T(S_1,S_2)$ admits a
%euclidean scalar product such that the different tensor products
%of the colors are forming an orthonormal basis. We call this basis
%the {\it color basis} of $T(S_1,S_2)$. Using the partition
%$P(S_1)\cup P(S_2)$ and the construction from section \ref{riwq}
%we get a linear function
%$$m(S_1,S_2):T(S_1,S_2)\rightarrow W_{S_1\cup S_2}\subseteq P_d.$$
%A natural orthonormal basis of $\hat{A}_d(S_1,S_2)$ is formed by
%the elementary matrices which correspond to colorings $O(S_1)\cup
%O(S_2)\rightarrow\mathcal{C}$. Since the elements of the color
%basis of $T(S_1,S_2)$ also correspond to colorings $O(S_1)\cup
%O(S_2)\rightarrow\mathcal{C}$ there is a natural isomorphism
%between the euclideal spaces $T(S_1,S_2)$ and
%$\hat{A}_d(S_1,S_2)$. By abusing the notation we identify these
%two spaces. Now there is unique map
%$$m:\hat{A}_d\rightarrow P_d$$ such that $m$ restricted to
%$\hat{A}_d(S_1,S_2)$ equals $m(S_1,S_2)$.

\subsection{The averaging operator}

We define the {\it averaging operator} $\xi:P_d\rightarrow I_d$ by
$$\xi(g)=\int_{\alpha\in O_d(\Bbb{R})}g^{\alpha}~d\nu$$
where $\nu$ is the normalized Haar measure on the orthogonal group
$O_d(\Bbb{R})$.
\begin{lemma}\label{commu}
The following diagram is commutative:

\begin{picture}(140,160)(10,10)

\put(80,140){$\hat{A}_{d,r}$}\put(102,142){\vector(1,0){100}}

\put(205,140){$P_d$}

\put(146,148){$m$}

\put(84,136){\vector(0,-1){45}}

\put(80,83){$\omega(A_{d,r})$}

\put(86,110){$\mathcal{P}_{d,r}$}

\put(116,85){\vector(1,0){86}}

\put(205,83){$I_d$}

\put(209,136){\vector(0,-1){43}}

\put(146,89){$m$}

\put(212,110){$\xi$}

\put(40,86){\vector(1,0){37}}

\put(20,83){$A_{d,r}$}

\put(32,23){\vector(3,1){170}}

\put(25,80){\vector(0,-1){50}}

\put(20,19){$\mathcal{Q}$}

\put(50,89){$\omega$}

\put(28,55){$\tau$}

\put(34,21){\vector(1,0){167}}

\put(205,19){$\Bbb{R}$}

\put(209,79){\vector(0,-1){50}}

\put(115,27){$f$}

\put(115,57){$t_d$}

\put(212,55){$\hat{f}$}
\end{picture}

\end{lemma}

\begin{proof}
Note that each map on the diagram is $\Bbb{R}$-linear and so it is
enough to check the commutativity for an appropriately chosen
generating system of the spaces.
\smallskip

First we prove that $\xi\circ m=m\circ \mathcal{P}_{d,m}$ by
checking it for the spaces $\hat{A}_d(S_1,S_2)$. Recall that
$\hat{A}_d(S_1,S_2)$ is identified with the euclidean space
$T=T(O(S_1)\cup O(S_2), P(S_1)\cup P(S_2))$ and that the
orthogonal group $O_d(\Bbb{R})$ is acting on $T$ by taking the
tensor product of the actions on $V$. This action commutes with
the map $m$ and so we have that
$$m\Bigl(\int_{\alpha\in O_d(\Bbb{R})}t^\alpha~d\nu\Bigr)=\xi(m(t))$$ for
all $t\in T$. Since the action of $O_d(\Bbb{R})$ preserves the
scalar product on $T$ one gets that $\int_{\alpha\in
O_d(\Bbb{R})}t^\alpha~d\nu$ is the orthogonal projection of $t$ to
the space of invariant elements. Therefore it suffices to prove
that $\omega(A_{d,r})$ is the space of invariants. This follows
from Weyl's first fundamental theorem as described in chapter
\ref{riwq}.
\smallskip

One gets $m\circ\omega=t_d\circ\tau$ by showing that
$$m(\omega(a(S_1,S_2,M)))=t_d(\tau(a(S_1,S_2,M)))$$ for all
triples $S_1,S_2,M$. This follows immediately from the
definitions.

The statement $\hat{f}\circ t_d=f$ is proved in corollary
\ref{clifting}.
\end{proof}

\begin{lemma}\label{posrep}
If $g\in P_d$ then for a sufficiently large natural number $r$
there is a symmetric positive semi-definite matrix $M$ in
$\hat{A}_{d,r}$ such that $m(M)=g^2$.
\end{lemma}

\begin{proof}
If $r$ is a large-enough natural number then
$$g=\sum_{S\in\mathcal{S},\mu(S)\leq r}g_S$$ where $g_S$ is an
element of $W_S$. Let us represent each $g_S$ by an element $t_S$
in the space $T(O(S),P(S))$ such that $m(t_S)=g_S$. This is
possible because $m$ is a surjective map to $W_S$. Setting
$$t=\sum_{S\in\mathcal{S},~\mu(S)\leq m}t_S\in\hat{B}_{d,m}$$ we
have that $m(t)=g$ and that $m(t\otimes t)=g^2$. On the other hand
$t\otimes t$ is represented as a rank $1$ positive semi-definite
matrix $M$ in $\hat{A}_{d,r}$.
\end{proof}

\begin{lemma}\label{negyzetposi}
If $g\in P_d$ then $\hat{f}(\xi(g^2))\geq 0$.
\end{lemma}

\begin{proof}
Using lemma \ref{posrep} we get that for a sufficiently large $r$
there is a symmetric positive semi-definite matrix
$M\in\hat{A}_{d,r}$ such that $m(M)=g^2$. From lemma \ref{projalg}
we obtain that $\mathcal{P}_{d,r}(M)=K^2$ where $K$ is a symmetric
matrix from $\omega(A_{d,r})$. Let $\bar{K}$ be a preimage of $K$
under the map $\omega$. We have that
$\omega(\bar{K}\bar{K}^T)=K^2$. By lemma \ref{posit} it follows
that $f(\tau(\bar{K}\bar{K}^T))\geq 0$. Lemma \ref{commu} implies
that
$$f(\tau(\bar{K}\bar{K}^T))=\hat{f}(m(\omega(\bar{K}\bar{K}^T)))=\hat{f}(m(K^2)).$$
It follows that $\hat{f}(m(\mathcal{P}_{d,r}(M)))\geq 0$. Using
lemma \ref{commu} again we obtain that $\hat{f}(\xi(m(M)))\geq 0$
which completes the proof.
\end{proof}

\subsection{Extension of $\hat{f}$ to $P_d$}

In this section we finish the proof of our main theorem by showing
that $\hat{f}:I_d\rightarrow\Bbb{R}$ extends to a homomorphism
$\bar{f}:P_d\rightarrow\Bbb{R}$. This is clearly enough because
the edge coloring model defined by
$$t(v)=\bar{f}(x_v),~~v\in\Bbb{N}^d$$ is a real valued edge coloring model
which represents the graph parameter $f$. We will need the
following well known consequence of the so-called
Positivestellensatz (see: \cite{RAG}).
\begin{theorem}\label{pst} Let
$g\in\Bbb{R}[x_1,x_2,\dots,x_n]$ be a polynomial such that it has
no root in $\Bbb{R}^n$. Then there exist polynomials
$p,f_1,f_2,\dots,f_h$ for some natural number $h$ such that
$$pg=1+f_1^2+f_2^2+\dots+f_h^2.$$
\end{theorem}

\bigskip

Let $P_{d,r}$ be the subring of $P_d$ which is generated by the
variables $\{x_v~|~h(v)\leq r,~v\in\Bbb{N}^d\}$. Since $P_{d,r}$
is the direct sum of the spaces $W_S$ where $S$ is a multiset of
$\{0,1,\dots,r\}$ we have that $P_{d,r}$ is invariant under the
action of $O_d(\Bbb{R})$. Lemma \ref{inrep} shows that
$I_{d,r}=I_d\cap P_{d,r}$ is the set of invariant polynomials in
$P_{d,r}$. Let $N_{d,r}$ be the kernel of the homomorphism
$\hat{f}:I_{d,r}\rightarrow\Bbb{R}$.

\begin{lemma}\label{firstext}
The homomorphism $\hat{f}:I_{d,r}\rightarrow\Bbb{R}$ extends to a
homomorphism $\bar{f}_r:P_{d,r}\rightarrow\Bbb{R}$.
\end{lemma}

\begin{proof}
First of all note that $P_{d,r}$ is a polynomial ring with
$t_r=\sum_{i=0}^r {{i+d-1}\choose{d-1}}$ variables. It suffices to
prove that there is a point $x$ in $\Bbb{R}^{t_r}$ which is a
common root for all the polynomials in $N_{d,r}$ because the
substitution of $x$ into polynomials from $P_d$ would yield a
homomorphism of the required form. Let $M$ be the ideal generated
by $N_{d,r}$ in $P_{d,r}$. Since $P_{d,r}$ is Noetherian we have
that there are finitely many polynomials $g_1,g_2,\dots,g_k\in
N_{d,r}$ which generate $M$ as an ideal. We prove by contradiction
that $g_1,g_1,\dots,g_k$ have a common root. Assume that it is not
true. Then $s=\sum_{i=1}^k g_i^2$ is a polynomial in $N_{d,r}$
which is positive everywhere in $\Bbb{R}^{t_r}$. Using theorem
\ref{pst} we get that there is a polynomial $p\in P_{d,r}$ such
that
$$ps=1+f_1^2+f_2^2+\dots+f_h^2$$ for some natural number $h$.
Applying the averaging operator $\xi$ for both sides we get that
$$\xi(p)s=1+\xi(f_1^2)+\xi(f_2^2)+\dots+\xi(f_h^2)$$ because $s$
is invariant under the action of $O_d(\Bbb{R})$. The left side is
an element of $N_{d,r}$ since $N_{d,r}$ is an ideal in $I_{d,r}$
and $\xi(p)$ is an element of $I_{d,r}$. This is a contradiction
because lemma \ref{negyzetposi} shows that
$$\hat{f}(1+\xi(f_1^2)+\xi(f_2)^2+\dots+\xi(f_h^2))\geq 1$$ which
means that the right side in not an element of $N_{d,r}$.
\end{proof}

\begin{lemma}
The map $\hat{f}:I_d\rightarrow\Bbb{R}$ extends to a homomorphism
$\bar{f}:P_d\rightarrow\Bbb{R}$.
\end{lemma}

\begin{proof}
Let $$g_s=\sum_{h(v)=s,~v\in\Bbb{N}^d}g_v^2.$$ It is easy to see
that $g_s=t_d(G_s)$ where $G_s$ is the graph with two nodes which
are connected by $s$ edges. This implies that $g_s$ is an element
of $I_d$. Let $\bar{f}_r$ be a map described by lemma
\ref{firstext}. If $r\geq s$ then $g_s$ is in $P_{d,r}$ and
$$f(G_s)=\hat{f}(g_s)=\sum_{h(v)=s,~v\in\Bbb{N}^d}\bar{f}_r(g_v)^2.$$
It follows that $|\bar{f}_r(g_v)|\leq\sqrt{f(G_s)}$ for all
$v\in\Bbb{N}^d$. Using these inequalities we have that there is an
infinite sequence $r_1<r_2<\dots$ of natural numbers such that
$\hat{f}_{r_i}(x_v)$ is convergent for all fixed vector
$v\in\Bbb{N}^d$. This means that $\hat{f}_{r_i}(p)$ is convergent
for all polynomial $p\in P_d$ and the limit is a homomorphism
which is an extension of $\hat{f}$ to $P_d$.
\end{proof}

\section*{Acknowledgement}
I am very indebted to L\'aszl\'o Lov\'asz for introducing me to
this beautiful subject, and for the many exciting discussions
without which this paper couldn't have been written. I am also
indebted to Michael H. Freedman and M\'aty\'as Domokos for their
kind help and suggestions.

\end{document}